\documentclass{aims}
\usepackage{amsmath}
  \usepackage{paralist}
  \usepackage{graphics} 
  \usepackage{epsfig} 
\usepackage{graphicx}  \usepackage{epstopdf}
 \usepackage[colorlinks=true]{hyperref}
\hypersetup{urlcolor=blue, citecolor=red}

\def\currentvolume{X}
 \def\currentissue{X}
  \def\currentyear{200X}
   \def\currentmonth{XX}
    \def\ppages{X--XX}
     \def\DOI{10.3934/xx.xx.xx.xx}

\newtheorem{theorem}{Theorem}[section]
\newtheorem{corollary}{Corollary}

\newtheorem{lemma}[theorem]{Lemma}
\newtheorem{proposition}{Proposition}

\theoremstyle{definition}
\newtheorem{definition}[theorem]{Definition}
\newtheorem{remark}{Remark}

\title[Near Force-Free Magnetic Fields, Blind ] 
      {On the Topological Characterization of Near Force-Free Magnetic Fields, and the Work of  Late-onset Visually-impaired Topologists }

\author[P. R. Kotiuga]{}

\subjclass{Primary: 58; Secondary: 78.}

 \email{prk@bu.edu}

\thanks{The author enjoys the multi-faceted support of his wonderful wife}

\begin{document}
\maketitle

\centerline{\scshape P. Robert Kotiuga }
\medskip
{\footnotesize
 \centerline{Boston University, ECE Dept.}
   \centerline{8 Saint Mary's Street}
   \centerline{ Boston, MA 02215, USA}
} 

\bigskip

 \centerline{(Communicated by the associate editor name)}

\begin{abstract}
The Giroux correspondence and the notion of a near force-free magnetic field  are used to topologically characterize  near force-free magnetic fields which describe a variety of physical processes, including plasma equilibrium. As a byproduct, the topological characterization of force-free magnetic fields associated with current-carrying links, as conjectured by Crager and Kotiuga, is shown to be necessary and conditions for sufficiency are given. Along the way a paradox is exposed: The seemingly unintuitive mathematical tools, often associated to higher dimensional topology, have their origins in three dimensional contexts but in the hands of late-onset visually impaired topologists. This paradox was previously exposed in the context of algorithms for the visualization of three-dimensional magnetic fields. For this reason, the paper concludes by developing connections between mathematics and cognitive science in this specific context.
\end{abstract}


\section{Introduction}

The primary mathematical result of this paper is theorem (\ref{InducedCorrespondence}). It exploits the Giroux correspondence to relate two characterizations of near force-free magnetic fields (in the absence of reconnection points). How the author came to this result over the course of decades is independent of the presentation given here. Over the last three decades, the finite element analysis of quasistatic  electromagnetic fields has moved from scalar problems in two dimensions to vector problems in three. Since two dimensional plots are understood cognitively in a manner quite different from our three dimensional experience, there were false starts.  Pitfalls have been avoided with great success by letting n-dimensional tools dominate algorithm development, but the results can threaten our {\em intuitive notions of visualization} - a term we dare not define with any precision.  Historically, it is surprising that the n-dimensional mathematical results used in this paper came into being by a very curious route. They originated in three dimensions, but in the hand of late-onset blind topologists. The results were then extended to higher dimensions where we are all equally blind. These surprising connections are explored in the second half of the paper and a hypothesis about this is made in the hope that it will help us reexamine the relationship between mathematics and our expectation of {\em visualization}. The paper is arranged as follows:

\tableofcontents

Concepts used are easily traced to original sources or textbooks; interdisciplinary connections are emphasized here. A unique approach is attempted in the hope of connecting with the intended audience.  First, there are no pictures in this paper\footnote{Unlike Lagrange's magnum opus, this is not because all problems have been reduced to calculation. Rather, the reader is encouraged to visualize the material in the manner he or she feels most comfortable. The sighted can find pictures on the web (and then question the extent to which such pictures help grasp the generality of the results!)}. Second, calculation and algorithms have taken a back seat in order to focus on three-dimensional topological issues\footnote{This is ironic since the subject is rooted in the use of algebraic topology in computational electromagnetics. Hopefully, this will be appreciated as one progresses through the paper.}. 

To make topological aspects which do not depend on the choice of Riemannian metric unambiguous, there is an interplay between the traditional vector analysis notation of electromagnetic theory and mathematical results stated in terms of the formalism of differential forms and manifold theory. Vector notation obscures distinctions between 1-forms and 2-forms in three dimensions, but reverting to the integral form of Maxwell's equations in $(3+1)$ dimensions resolves any ambiguity. 

\section{The Lorentz force law and twisted magnetic fields}

The Lorentz force,
\begin{equation}
{\bf F}={\bf J}\times{\bf B},
\end{equation}
relates the magnetic force density ${\bf F}$, the current density ${\bf J}$ and the magnetic flux density ${\bf B}$.

\begin{definition} \label{force-free def}
 Consider the linear combination
$a{\bf J}+b{\bf B}$
\begin{itemize}
\item A {\em force-free magnetic field} is a magnetic field whose Lorentz force vanishes. That is, the above linear combination of ${\bf J}$ and ${\bf B}$  vanishes identically for some scalar function coefficients $a$ and $b$. 
\item A {\em linear force-free magnetic field} is a magnetic field where the above linear combination of ${\bf J}$ and ${\bf B}$ vanishes identically with coefficients $a$ and $b$ constant. 
\end{itemize}
\end{definition}

\begin{remark}
Let $\mu$ and $\bf H$, defined implicitly below, be the permeability and the magnetic field intensity respectively.
Nontrivial linear force-free magnetic fields are a popular notion in the magnetohydrodynamics literature. They can result from any eigenfield of the curl operator and the quasistatic version Ampere's law, $\rm curl\bf H =\bf J$,
\begin{equation}
\bf J\times\bf B=(\rm curl\bf  H)\times\mu\bf H=\lambda\mu(\bf  H\times\bf  H)=0
\end{equation}
We will see that the linear force-free magnetic field problem has a nice variational formulation via the the eigenvalue problem for a self-adjoint curl operator.
\end{remark}

The bridge from force-free magnetic fields (or the more restrictive notion of linear force-free fields) to the global topological aspects is provided by what the author calls ``trigonometry with a twist". Specifically, recall the most basic of trig identities:
\begin{equation}
 1=\sin^2\theta+\cos^2\theta
\end{equation}
In the context of magnetic fields interpret this as:
\begin{equation} \label{twisted trig}
|\bf J|^2|\bf B|^2=|\bf J\times\bf B|^2+(\bf J\cdot\bf B)^2
\end{equation}
The ``twist" in this equation will be seen to follow from the violation of the ``Frobenius Integrability condition".

\begin{remark} \label{Integrabilitiy}
 Under mild hypotheses about differentiability, it is a simple exercise in partial differentiation to verify the following statement:
``If ${\bf F} = g \rm grad f$ for some functions $f$ and $g$, then $\bf F \cdot \rm curl\bf F = 0$". However, the converse statement, 
``If $\bf F\cdot\rm curl\bf F = 0$ in $R$, then ${\bf F} = g \rm grad f$ for some $f$ and $g$",
 is not true in general. Furthermore, because the constraint is quadratic, there are no general global conditions on the topology of the domain $R$, which ensure that the converse holds. The general conditions depend on the topology of the field lines of $\bf F$. In fact, one can show that for any convex region, one can find an $\bf F$ such that $\bf F\cdot\rm curl\bf F=0$ but $\bf F$ cannot be represented as $g \rm grad f$\footnote{This is in contrast to the special case of defining a magnetic scalar potential, where the constraint becomes linear and additional hypotheses about (co)homology groups provided necessary and sufficient conditions for the converse to be true.}. It turns out that if $\bf F$ is nonzero in a neighborhood of some point, then there is a smaller neighborhood in which the converse is true. This is formalized by the Frobenius integrability theorem...
\end{remark}

\begin{lemma}
 In terms of differential forms, the {\bf Frobenius integrability condition} is
$\omega \wedge d\omega = 0$, where $\omega$ can be a k-form on an n-dimensional manifold, and $d$ is the exterior derivative. Hence, the condition is independent of any metric notion. The {\bf Frobenius integrability theorem} states that when the Frobenius integrability condition is satisfied in a an open set $U$ containing a point $p$, then there is a subset $V$ of $U$ containing $p$ such that $d\omega =\eta \wedge \omega$ for some closed 1-form $\eta$. 
\end{lemma}

For $k=1$, the closed form is locally exact, so a coordinate calculation involving logarithmic derivatives and a translation into classical vector notation yields...
\begin{corollary}
For 1-forms in 3-d Euclidean space and adopting classical vector notation, the Frobenius integrability condition is $\bf F\cdot\rm curl\bf F = 0$.  The Frobenius integrability theorem implies that when the Frobenius integrability condition is satisfied in a an open set $U$ containing a point $p$, {\bf and} $\bf F$ is nonvanishing on $U$, then there is a subset $V$ of $U$ containing $p$ such that ${\bf F} = g \rm grad f$ for some functions $f$ and $g$.   
\end{corollary}

\begin{remark} \label{rescaled equipotentials}
From this corollary we see that if the Frobenius condition is satisfied locally by a nonvanishing vector field $\bf F$, then $\bf F$ can still be expressed in terms of a scalar potential $f$ and one still has the notion of  equipotentials, but the equipotentials are rescaled by an {\em integrating factor} $g$. The hypothesis that $\bf F$ is locally nonvanishing, ensures that these rescaled equipotentials are nonsingular. As such, these rescaled equipotentials are called the {\em leaves} of a {\em foliation}. That is, a foliation is the result of considering a region of space to be a stack of leaves given by the rescaled equipotentials of $\bf F$. In our case, the level sets of the function $f$ are locally, the
leaves of the foliation. Note that on sufficiently small neighborhoods these leaves are compact, but they may not be compact globally since the functions $f$ and $g$ are only defined on coordinate patches yet the leaves are well-defined globally\footnote{ Note that neighborhoods where $\bf F$ vanishes need to be excised and in the process, further global considerations can arise.}.
\end{remark}

\begin{remark} \label{propeller}
When $\bf F\cdot\rm curl\bf F \ne 0$, the Frobenius condition is violated, the notion of rescaled equipotentials are no longer applicable, there can be no foliation, but a topological interpretation of this situation is crucial to further developments. If $\bf F\cdot\rm curl\bf F \ne 0$ in a neighborhood of a point $p$ then it cannot switch sign in this neighborhood. If one tries to define equipotentials in infinitesimal neighborhoods of points in this neighborhood then they are nonintegrable and line up like the fins of a propeller. Hence $\bf F\cdot\rm curl\bf F$ measures the {\em twisting} of the vector field. This notion of twisting will be key to defining {\em near force-free} magnetic fields and their relation to the Giroux correspondence.
\end{remark}

We are now positioned to see how twisting, in a well-defined topological sense,  is related to quantities familiar from Maxwell's electrodynamics. Consider for simplicity Ohmic conductors ($\bf J = \sigma \bf E$),  an absence of magnetizable media ($\bf B = \mu \bf H$) and negligible displacement currents. Here $\sigma$ is the conductivity and $\bf E$ the electric field intensity. Equation (\ref{twisted trig}) then becomes:
\begin{equation} \label{twisted trig n phys}
\sigma\mu(\bf J\cdot\bf E)(\bf B\cdot\bf H)=|\bf J\times\bf B|^2+\mu^2(\bf H\cdot\rm curl\bf H)^2
\end{equation}
The l.h.s. of eq(\ref{twisted trig n phys}), involves the densities of Joule heating $(\bf J\cdot\bf E)$, and magnetic energy $(\bf B\cdot\bf H)$, while the r.h.s. involves the magnitude of the Lorentz force density $(\bf J\times\bf B)$ and a local measure of topological twisting $(\bf H\cdot\rm curl\bf H)$. ``Trigonometry with a twist" says that in regions where the magnetic field and heating are nonzero, one must have either a nontrivial Lorentz force or twisted field lines. In this paper we exploit this local identity and the notion of a confoliation, to  globally characterize the topology of {\em near force-free} magnetic fields. The notion of {\em current helicity} and the eignvalue problem for linear force-free fields will play a key role.

\begin{definition}
If we identify the magnetic field intensity as a 1-form, $\omega=\bf H\cdot dr$, the {\em current helicity functional} associated with the magnetic field on a region $R$ is
\begin{equation}
\mathcal{H}(\rm curl\bf H)=
\frac{1}{2}\int_{{I\!\! R}^3}\bf H\cdot\rm curl\bf H dV=
\frac{1}{2}\int_{{I\!\! R}^3}\omega\wedge d \omega 
\end{equation}
\end{definition}

\begin{remark} 
The fact that the current helicity functional, when expressed in terms of differential forms, depends only on the exterior derivative and the wedge product but not on the Hodge star, leads to important invariance properties:
\begin{enumerate}
\item This integral is diffeomorphism invariant. However, in the physical interpretation of the current helicity, it is the current which has to be identified with the 2-form  $d\omega$. This is done via a Lie derivative, by contracting a vector field with the volume form. As a result, as a function of the current density, {\em the integral is invariant under volume preserving diffeomorphisms}, a Lie group which we denote by $Diff_0(R)$\footnote{We denote the tangent space to $Diff_0(R)$ by $T(Diff_0(R))$. This ``Lie algebra of infinitesimal volume preserving diffeomorphisms" is precisely the space of divergence-zero vector fields.}. What is striking is that as a Lie algebra, $T(Diff_0(R))$ has an ``invariant bilinear form" and it is  precisely the bilinear form associated with the helicity functional. This bilinear form has a conceptually satisfying topological interpretation as Arnold's "mean asymptotic linking number"\footnote{To make a precise statement, boundary conditions need to be specified along with additional hypotheses in the case of $R$ not being simply connected.}.
\item The integral is {\em gauge invariant} in the sense that if one adds an exact 1-form to $\eta$, an ``integration by parts"  shows that the integral changes by a boundary term which vanishes under suitable Lagrangian boundary conditions. Since this is independent of any Riemannian metric, it has interesting consequences in the context of Hodge decompositions. 
\end{enumerate}
\end{remark}

We will now see how the current helicity enables to construct linear force-free fields arise as eigenfields of the $\rm curl$ operator by means of a variational formulation. Relative to the Hodge inner product, 
\begin{equation}
\langle{\bf E},{\bf F}\rangle=\int_R {\bf E}\cdot{\bf F} dV,
\end{equation}
an ``integration by parts" formula shows that curl is formally self-adjoint\footnote{The construction of self-adjoint curl operators is the subject of another paper (\cite{HKT})}:
\begin{equation}
\langle{\bf D},\rm cur{l\bf C}\rangle=\langle\rm curl{\bf D},{\bf C}\rangle
  +\int_{\partial R}({\bf D}\times{\bf C})\cdot\widehat{\bf n}dS
\end{equation}
If $R$ is a compact manifold with boundary, the boundary term is a symplectic form. If we impose ``Lagrangian boundary conditions", the boundary term vanishes and the curl operator becomes a self-adjoint Fredholm operator in some Sobolev-like space where the gradients have been ``quotiented out". We then have a countable number of finite energy divergence-zero eigenfields satisfying an equation of the form:
\begin{equation}
\rm curl\bf H = \lambda\bf H
\end{equation} 

\begin{remark} \label{HelicityEVP}
In terms of the current helicity functional, the linear force-free magnetic field problem now has variational formulation. That is, to find the stationary points of $\mathcal{H}(\rm curl\bf H)$ subject to Lagrangian boundary conditions and the normalization $(\langle{\bf H},{\bf H}\rangle=1)$ imposed via a Lagrange multiplier. The eigenfield corresponding to the nonzero eigenvalue having the smallest absolute value can be interpreted as either the magnetic field which has the most twisting for a given amount of energy, or the magnetic field with the least energy for a given amount of twisting.
\end{remark}
\begin{remark}
What is unexpectedly elegant is that when the helicity functional is restricted to a finite domain, the diffeomorphism-invariant Lagrangian boundary conditions which render the curl operator self-adjoint depend only on the Lagrangian subspaces of the first homology group of the boundary. See (\cite{HKT}) for details. Even more rewarding is that all these properties have a discrete counterpart under Whitney-form finite element discretization.
\end{remark}
\begin{remark} \label{BigEVPplans}
By restricting the curl operator to subdomains and solving finding the eigenfunctions of the curl operator on these subdomains, one can construct a global force-free magnetic fields out of piecewise linear force-free fields. This technique will play a key role in reconciling the Giroux correspondence with energy constraints. Specifically, we plan to apply these ideas to the {\em binding of an open book}. 
\end{remark}
 Before developing a topological characterization of force-free magnetic fields in terms of ``open books", we now relate foliations to notions from electromagnetics.

\subsection{Cuts for magnetic scalar potentials and preferred foliations }

We now set the stage for confoliations and the Giroux correspondence by focusing on particular foliations associated with maps into circles. These maps have an interpretation in terms of (multi-valued) magnetic scalar potentials which arise in the calculation of ``cuts" for magnetic scalar potentials. It turns out that the set of ``regular values" of these maps yield foliations with ``compact leaves", and if the map has no critical points, the entire region is foliated by compact leaves.

To get started, we need a formulation of de Rham's theorems\footnote {De Rham's theorems were originally formulated as two theorems; one dealing with the injectivity of a map and the other dealing with surjectivity. Since the coefficient ring is a field in the lemma, de Rham's theorems are an isomorphism between a vector space and its dual. We have not specified the space of (co)chains because for any reasonable choice, the Eilenberg-Steenrod axioms for a (co)homology theory will bail us out if we get into a trouble.},\footnote{We also avoid stating the corresponding de Rham isomorphisms for relative (co)homology groups of $M$ modulo boundary since the mathematical development doesn't require them except in the proof of lemma(\ref{Lefschetz}). They can arise by considering differential forms with compact supports on an open manifold with compact closure. In this paper, magnetic flux linking current paths are examples of relative periods.}.
\begin{lemma}
Let $M$ be a n-dimensional compact orientable manifold with (possibly empty) boundary and ${I\!\! F}$ ($ ={I\!\! R} $ or $  {I\!\!\!\! C}$), a field, then
\begin{equation}
H^{k}(M;{I\!\! F}) \simeq H_{k}(M;{I\!\! F})
\end{equation}
Denote by $[c]$ and $[\omega]$ the (co)homology classes associated with the (co)cycles $c$ and $\omega$ respectively. The isomorphism is given by the period map,
\begin{equation}
 [\omega, c] \equiv \int_{c}{\omega} : H^{k}(M;{I\!\! F}) \times H_{k}(M;{I\!\! F}) \longrightarrow {I\!\! F}
\end{equation}
considered as a well-defined nondegenerate bilinear pairing on (co)homology.
\end{lemma}

For pedestrians: The hypothesis of mixed second partial derivatives commuting, implies the following statement:
``If  $\bf F = \rm grad f$ for some differentiable function $f$ on $R$, then $\rm curl\bf {F=0}$ in $R$".
In the following corollary, modeled on remark(\ref{Integrabilitiy}), the conditions for the existence of scalar potentials is stated in terms of the converse.

\begin{corollary} 
Let $R$ be any compact orientable three-dimensional manifold with boundary. The statement ``If $\rm curl\bf { F = 0}$ in $R$, then $\bf F = \rm grad f$ for some $f$", is true
\begin{itemize}
\item  Whenever $R$ is simply connected and more generally if the first homology group of the region $H_{1}(R;{I\!\! R})$, is trivial. 
\item If $H_{1}(R;{I\!\! R})$ is nontrivial, but all the periods of $ \omega = \bf {F\cdot dr}$ vanish. 
That is, the cohomology class of $\omega$ is zero.
\end{itemize} 
\end{corollary}

\begin{remark}
When the coefficient ring is a principal ideal domain (e.g. the integers ${Z\!\!\! Z}$), we no longer have a vector space, but a module and there might be torsion. In computational electromagnetics and in the computation of cuts for magnetic scalar potentials one needs integer (co)homology groups and so we refer to a standard lemma:  
\end{remark}
\begin{lemma}
For $R$, any compact orientable three-dimensional manifold with boundary embedded in ${I\!\! R}^3$, all of the torsion subgroups  of  the (co)homology groups of $R$, or $R$ modulo its boundary, when computed with coefficients in ${Z\!\!\! Z}$, vanish.
\end{lemma}
Minimal effort then produces a reformulation which is very useful for this paper: 
\begin{corollary} \label{tors free and circle maps}
For $R$, any compact orientable three-dimensional manifold with boundary embedded in ${I\!\! R}^3$, for any (co)homology group of $R$, or $R$ modulo its boundary, let $\beta$ be the dimension of the resulting vector space when the underlying coefficient group is ${I\!\! R}$. Then if  the underlying coefficient group is ${Z\!\!\! Z}$, the resulting abelian group is free, is isomorphic to ${Z\!\!\! Z}^{\beta}$, and forms a lattice in the vector space.

Suppose further that $R$ is connected with distinguished point $p_0$. Then, for any closed 1-form $\omega$ (i.e. 1-cocycle) with integer periods, there is a well defined map from $R$ to the unit circle in the complex plane, $ f:R\longrightarrow S^{1}$, given by
\begin{equation}
f(p)= exp({2\pi i}\int_{p_0}^{p} \omega)
\end{equation}
and $[M, S^{1}]$, the homotopy classes of such maps, form a group isomorphic to  $H^{1}(R;{Z\!\!\! Z})$.
\end{corollary}

\begin{remark} \label{curl zero}
Corollary(\ref{tors free and circle maps}) yields special foliations. Scalar potentials arise when the Frobenius condition is satisfied trivially because $\rm curl\bf F = 0$ and the periods of $\bf {F\cdot dr}$ vanish. In this case, the integrating factor of remark(\ref{Integrabilitiy}) disappears (i.e. $g=1$), and equipotentials are compact. If the periods of $\bf {F\cdot dr}$ do not vanish, then there is no globally defined scalar potential. However, there can be a singular foliation whose leaves are compact if and only if the periods of $\bf {F\cdot dr}$ are ``commensurable". That is, the ratio of any two periods is a rational number. In this case, we can rescale all of the periods, obtain an integer cohomology class and using the map constructed in the corollary, define the leaves of the singular foliation to be the inverse images of points on the circle. 
This observation lets us define maps into circles whose regular values form a foliation with compact leaves\footnote{These maps, which appeared in computational electromagnetics over a quarter century ago in the context of {\em cuts} for magnetic scalar potentials (\cite{cuts1}), are now crucial to our use of the Giroux correspondence via the notion of an open book decomposition.}.
\end{remark}
To understand why ``cuts" for magnetic scalar potentials are related to special foliations, we will need to introduce Lefschetz duality.
\begin{lemma} \label{Lefschetz}
Let $M$ be a n-dimensional compact connected orientable manifold with boundary $\partial M$, and ${Z\!\!\! Z}$, the ring of integers. Then $M$ has a fundamental homology class which generates $H_{n}(M, \partial M;{Z\!\!\! Z})$ and there is a natural isomorphisms
\begin{equation}
H^{n-k}(M, \partial M;{Z\!\!\! Z}) \simeq H_{k}(M;{Z\!\!\! Z}),
\end{equation}
\begin{equation}
H_{n-k}(M, \partial M;{Z\!\!\! Z}) \simeq H^{k}(M;{Z\!\!\! Z}).
\end{equation}
These isomorphisms arise from two distinct realizations of the modules dual to $H^{k}(M;{Z\!\!\! Z})$ and $H^{n-k}(M, \partial M;{Z\!\!\! Z})$ respectively. One realization comes from the bilinear pairing ``induced by integration" and the other, the cup product evaluated on the fundamental class.  
\begin{equation}
 [ H_{n}(M, \partial M;{Z\!\!\! Z}), H^{n-k}(M, \partial M;{Z\!\!\! Z}) \cup H^{k}(M;{Z\!\!\! Z}) ] \longrightarrow {Z\!\!\! Z}
\end{equation}
\end{lemma}
\begin{remark}
Lefschetz duality,  being stated in terms of integer coefficients without mentioning torsion, is a {\em natural} isomorphism. It is valid for coefficients in any abelian group or any field; if the field is ${I\!\! R}$, then it applies to de Rham cohomology.
\end{remark}

\begin{corollary}
Let $R$ and $f(p)$ be as in corollary(\ref{tors free and circle maps}), $p$ a regular value of $f$, then 
\begin{equation}
[M, S^{1}] \simeq H^{1}(R;{Z\!\!\! Z}) \simeq H_{2}(R, \partial R;{Z\!\!\! Z}). 
\end{equation}
 The isomorphisms are given by corollary(\ref{tors free and circle maps}) on the l.h.s. and Lefschetz duality on the r.h.s.. Here, $f(p)^{-1}$ is an orientable embedded manifold whose relative homology class in $H_{2}(R, \partial R;{Z\!\!\! Z})$ is the Lefschetz dual in $H^{1}(R;{Z\!\!\! Z})$, to the  class of \begin{equation}
\omega =\frac{1}{2\pi i} d(ln(f)).
\end{equation} 
\end{corollary}
Thus corollary(\ref{tors free and circle maps}) and Lefschetz duality, Lemma(\ref{Lefschetz}), allow us to give a formal definition of a cut for a magnetic scalar potential as the inverse image of a regular value of a map representing a cohomology class.
\begin{remark}
If we think in terms of the exterior of a tubular neighborhood of a current carrying knot, there are three very useful ways to think of a cut:
\begin{itemize}
\item It enables one to use a magnetic scalar potential exterior to the cut. the jump in the potential as one crosses the cut is equal to the current in the wire.
\item A surface used to compute the magnetic flux linking the current carrying knot (assuming one does not have access to a vector potential). 
\item A Seifert surface for the knot!
\end{itemize}
\end{remark}
\begin{remark} \label{Pontryagin}
This formulation of the cuts problem leads to a first encounter with a blind topologist. Poincar\'e assumed (incorrectly) that a cycle representing a homology class could be realized as a manifold. Pontryagin, in a 1931 paper asked and partially answered this question in the category of differentiable manifolds. He did this in knot theory where he proved the existence of an orientable, embedded compact surface whose boundary is a given knot (i.e. a Seifert surface) years before Seifert's 1934 paper. His approach was far more conceptual, capable of generalization to higher dimensions, and points to the ``Pontryagin-Thom construction", yet Seifert's proof is the one found in knot theory books! It is intuitive, and easily presented on the page to a sighted audience, yet Pontryagin's approach constructs a family of cuts and the spectrum of {\em Oriented Bordism} (a generalized cohomology theory), and lays the groundwork for Rene Thom's Fields medal in the 1950s.
\end{remark}

Cuts are highly nonunique since one is free to choose any map in a homotopy class, and then any regular value of the map. Regardless of which map we chose, the inverse images of regular values form a foliation. Getting ambitious, we can seek ``nice maps" by putting an energy functional on the space of maps. For example, consider finding a solution to the variational problem of minimizing
\begin{equation} 
        F(f)=\int_{R}\rm grad \bar{f}\cdot\rm grad f\;dV,
\end{equation}
subject to $\bar{f}f=1$ in $R$, the constraint that constrains the image to $S^{1}$. For the $j$th cut, \(1\leq j\leq\beta_{1}(R)\), we also fix the homotopy class of the cut
\begin{equation}
        \frac{1}{2\pi i}\oint_{c_k}\rm grad(\ln f)\cdot dl=\delta_{jk}.
\end{equation}
 Here $1\leq k\leq\beta_{1}(R)$, and the  $c_{k}$, form a collection of curves representing a basis for $H_{1}(R;{Z\!\!\! Z})$. The solution to the above problem defines a harmonic map. Indeed, the angle of the image is a harmonic function on $R$ and readily interpretable in terms of magnetic fields\footnote{Algorithms for computing cuts based on the finite element method work directly with the Dirichlet integral and the(multi-valued) angle on the covering space of the circle(\cite{cuts4},\cite{gross}).}. Cuts computed in this fashion enjoy all the properties of level sets of elliptic equations, but their properties are hard to articulate\footnote{We will return to the problem of ``nice cuts" in section 3.}.

\subsection{ Confoliations and near force-free magnetic fields}

One way of getting a feel for contact structures and confoliations is to start locally with coordinate expressions for Monte and Clebsch potentials.

\begin{remark} \label{Monge}
An arbitrary vector field $\bf F$ does not satisfy the Frobenius integrability condition. Suppose however we found a $\bf F'$ that differs from $\bf F$ by a gradient, $\bf F' = \bf F - \rm grad h$,  and we require $\bf F'$ to satisfy the Frobenius condition, $\bf F'\cdot\rm curl\bf F'=0$. The function h would then  have to be a  solution to a linear first order p.d.e. whose coefficients depend on $\bf F$. Once such an h is found, $\bf  F'$ can be expressed in terms of two scalar functions and $\bf F$ can be written as:
\begin{equation}
{\bf F} = g \rm grad f + \rm grad h.
\end{equation}
Such a triple of functions is known as  {\em Monge potentials}. They exist locally under the same hypotheses as the Frobenius integrability theorem and are highly non-unique.
\end{remark}

\begin{remark} \label{Clebsch}
Monge potentials ensure that any vector field with zero divergence can be expressed
locally as ${\rm curl}{\bf F} = {\rm grad} g \times \rm grad f$. Here, $g$ and $f$ are known as {\em Clebsch potentials}. This leads to a new type of foliation; $\rm curl\bf F$ foliates, $R$, the region of definition and the ``leaves'' are the curves formed locally by intersecting level sets of $f$ and $g$.
In this way, Clebsch potentials parameterize flux tubes in solenoidal fields.  This picture is intimately tied to Faraday's ``tubes and slices''.
\end{remark}

\begin{remark} \label{triple product}
Locally Monge potentials give an interesting triple product:
\begin{equation} 
{\bf F}\cdot{\rm curl}{\bf F} = 
\left( {\rm grad} g \times \rm grad f\right) \cdot \rm grad h
\end{equation}
Using differential forms, if $\omega=\bf {F\cdot dr}$, this becomes manifestly metric invariant:
\begin{equation}
\omega \wedge d\omega = dg \wedge df \wedge dh
\end{equation}
\end{remark}

\begin{definition} \label{contact structure n foliation}
Let $\omega$ be a 1-form on $R$, a compact orientable 3-d manifold. Express the 3-form of the Frobenius condition as a multiple of the volume form:
\begin{equation}
\omega \wedge d\omega = mdV
\end{equation}
\begin{itemize}
\item $\omega$ is a contact structure if $m > 0$ everywhere, (or $m<0$ everywhere),
\item $\omega$ is a confoliation if $m \ge 0$ everywhere, (or $m\le 0$ everywhere). 
\end{itemize}
\end{definition}

\begin{remark}
The intuition of behind a contact structure was given in remark(\ref{propeller}) while the intuition of behind a foliation was given in remarks(\ref{rescaled equipotentials}) and (\ref{curl zero}). Naively, a confoliation is part contact structure and part foliation. A powerful perspective is the following: {\em Confoliations are the boundary of the space of contact structures}.
\end{remark}

Let us denote by $Supp(\bf F)$, the support of a vector field $\bf F$.

\begin{definition} \label{near force-free def}
 Let $V= Supp({\bf B}) \cap Supp({\bf J})$. The magnetic field
$\bf B$ is a {\em near force-free field} if on every connected component of $V$,
\begin{equation} 
|\bf J|^2|\bf B|^2 > |\bf J\times\bf B|.
\end{equation}
\end{definition}

Comparing definitions (\ref{force-free def}) and (\ref{near force-free def}) we have:
\begin{lemma} \label{inclusions}
The various notions of constrained magnetic fields are related by:\\
\center{{\bf  Linear force-free $\subset$ force-free $\subset$ near force-free}}
\end{lemma} 

By definition(\ref{near force-free def}) and equation(\ref{twisted trig}) we have
\begin{lemma} \label{near f-f n topology}
Let $\bf B$ be a near force-free field and suppose that 
\begin{itemize} 
\item  displacement currents are negligible so that $\rm curl\bf H =\bf J$
\item $Supp({\bf B}) = Supp({\bf H})$ so that $V= Supp({\bf H}) \cap Supp({\bf J})$.
\end{itemize}
then on any connected component of $V$, $\bf H\cdot\rm curl\bf H$ is nonzero and can't change sign.
\end{lemma}

From remark(\ref{curl zero}), definition (\ref{contact structure n foliation}) and lemma(\ref{near f-f n topology}),we now have
\begin{lemma} \label{near f-f n confoliations}
With the same hypotheses as lemma (\ref{near f-f n topology}), and if the sign of the current helicity density is same on each connected component of $V$, then
\begin{itemize}
\item $\bf H$ defines a singular foliation in the complement of $Supp({\bf J})$, 
\item  {\bf $\bf H$ defines a contact structure on $V$}
\end{itemize}
Furthermore, if ${\bf H}$ is nonzero on the complement of $Supp({\bf J})$, then
\begin{itemize}
\item $\bf H$ defines a foliation in the complement of $Supp({\bf J})$, 
\item  {\bf $\bf H$ defines a confoliation on ${I\!\! R}^3$}.
 \end{itemize}
\end{lemma}

\begin{remark}
 Lemmas (\ref{near f-f n confoliations}) and (\ref{inclusions}) make clear that confoliations and the mathematical results concerning them, are applicable to the topological characterization of near force-free magnetic fields and other more restrictive force constraints\footnote{These include, plasmas equilibria in either controlled environments or in space (here force constraints are implicit in the notion of equilibrium), superconducting magnet design (where the {\em normal temperature} is dependent on the magnetic field transverse to the current flow), and other applications. What may be surprising is that the additional hypotheses in lemma(\ref{near f-f n confoliations}) relate to ``reconnection points" in plasma physics.}. 
\end{remark}

 From a topological point of view, what remains unsatisfying is the lack of characterization of the leaves of the foliation in the complement of $Supp({\bf J})$. We also don't have a  simple characterization of the global topology of the magnetic field.
Open book decompositions and the Giroux correspondence will help resolve this.

\subsection{A toy inverse problems for near force-free magnetic fields}

Let's propose a problem which is rhetorical in the sense that our goal is not to try proving existence of solutions as a first step. Rather, we will use it to introduce open book decompositions and the Giroux correspondence in the context of near force-free magnetic fields. It is reasonable to expect the inverse problem to be formulated as a  ``free boundary value problems" for a system of PDEs. The toy problem is introduced here is an attempt to abstract the topological essence of equilibrium configurations by a reformulation in terms of knot and link theory, a field with well-studied and readily computable algebraic invariants.

\begin{remark}  
Although introducing knots and links has intuitive appeal, later grief can only be avoided with proper definitions and well-defined concepts. Under stereographic projection, there is a {\em one point compactification} of ${I\!\! R}^3$ where $({I\!\! R}^3) \cup \{ \infty \} $ is identified with ${S^{3}}$, the unit sphere in ${I\!\! R}^4$. Hence, if noncompactness  is an issue in subsequent developments,  identify ${I\!\! R}^3$ with ${S^{3}}$. An  {\em isotopy} is a homotopy through embeddings. Knots are embeddings of $S^1$ into ${S^{3}}$ modulo an equivalence relation which we would like to be {\em ambient isotopy}. However,  this equivalence relation would make knot theory trivial. Hence, we define ``tame" knots and links: 
\end{remark}

\begin{definition}
An {\em n-component link} is an embedding of a disjoint union of n circles into ${I\!\! R}^3$ or ${S^{3}}$; $n=1$ gives a {\em knot}.
A {\em tame n-component link} is a link whose image is ambient isotopic to $n$ piecewise-linear curves; $n=1$ gives a {\em tame knot}.
\end{definition}

\begin{remark}  
 In subsequent developments, depending on context, by ``knots and links"  we will mean either ``tame knots and links" or their ambient isotopy classes.The definition of a tame knot also ensures the existence of a {\em tubular neighborhood}. Again, depending on context, when may refer to knots and  links, and their tubular neighborhoods interchangeably. This identification allows us to talk about ``current carrying knots and links" without a discussion about whether this is possible with finite energy.
\end{remark}

\begin{center}
 {\Large The Toy Inverse Problem for Current-carrying Links:}\\
{\bf Arrange a current-carrying link so the Lorentz force vanishes on it}. 
\end{center}

 Recall that a {\em closed braid}\cite{BurdeZeischang} is a link that can be represented in cylindrical coordinates as a curve whose tangent vector has a positive projection on the unit vector in the circumferential direction. That is, a closed braid generalizes the notion of a wound coil.  To expose the practical implications of this toy problem, consider Alexander's Theorem\cite{BurdeZeischang}:
\begin{lemma}
Any  link can be represented by a closed braid.
\end{lemma}

\begin{remark}  
By Alexander's theorem, we can isotope any knot into a closed  braid and so we ask: What ``knot types'' produce nice force-free coils?  Crager and Kotiuga conjectured that knot types of such optimal coils have equal Alexander and Thurston Norms.  See\cite{crager} for some pictures suggesting that torus knots both have the desired topological properties and a close connection to the field lines associated to the eigenfunctions of the curl operator on the torus. As we shall soon see ``open books'' are solutions to our toy problem, and by the theory of the Alexander polynomial\cite{BurdeZeischang}, 
 have the conjectured topological property. 
\end{remark}

Remark (\ref{curl zero}) shows that for a current-carrying link with ``commensurable" currents, the multivalued magnetic scalar potential in the exterior of the link has compact equipotentials. When the link is identified with $V$ in  lemmas(\ref{near f-f n topology}, \ref{near f-f n confoliations}), we have the conditions under which we have a confoliation and a near force-free magnetic field as a solution to our toy problem. Explicitly, we have

\begin{lemma} \label{ToyProbSolution}
Suppose the multivalued scalar potential has ``commensurable"  periods in the  link's complement, no critical points (so all equipotentials are isotopic), and the current forms a near force-free magnetic field in some tubular neighborhood of the link. Then the toy problem has a solution in the form of a confoliation.
\end{lemma}

\begin{remark} \label{PerfectEqNoReconnection}
In the link complement, lemma(\ref{ToyProbSolution}) yields solutions to our toy problem as circle-valued Morse functions\cite{Goda}, which are perfect (or perfect Morse-Novikov functions). In the language of space plasmas, near force-free magnetic fields are thus associated with the absence of {\em reconnection points}.
\end{remark} 

\subsection{Open-book decompositions and the Giroux correspondence } We follow  Etnyre\cite{Etnyre} for basic definitions and proofs pertaining to open book decompositions and the Giroux correspondence.

\begin{definition} An {\em open book decomposition} of a compact orientable manifold $M$ is a pair $(B, \pi)$ where:
\begin{itemize}
\item $B$ is an oriented link in M called the binding of the open book\footnote{ One might be tempted to call $B$ the ``spine" of the open book, but in the context of three dimensional manifolds, the ``spine of a manifold" has a specific technical meaning and so ``binding" is used in order to avoid confusion. However, since the foliation also defines a fiber bundle, one must remember that $B$ is not the base of the fiber bundle, but the binding of the open book.}.
\item and $\pi$ is a fibration from the complement of $B$ to $S^1$
  such that $\pi^{-1}(p)$ is the interior of a compact surface
  $\Sigma_p$ in $M$ and $\partial\Sigma_p=B$ for all $p$ in $S^1$. For any
  $p$ the surface $\Sigma_p$ is called a ``page'' of the open book. 
\end{itemize}
\end{definition}

\begin{corollary} \label{ReformulatedToyProbSolution}
of lemma(\ref{ToyProbSolution}).
Under stereographic projection of Euclidean space, solutions to our toy inverse problem which are free of ``reconnection points" are topologically characterized by open book decompositions of $S^3$ subject to positive twists on their bindings induced by the existence of a contact structure on a tubular neighborhood of the binding.
\end{corollary}

Two results give intuition for open book decompositions as special confoliations and their relationship to the space of contact structures. First, the connected components of the space of contact structures are contractible and confoliations which are not contact structures lie in the boundary\cite{Giroux2002}. Second, by a result of Thurston and Winkelnkemper\cite{ThurstWink}, every open book {\em supports} a contact structure (and so are distinguished confoliations). That is, an open book can deformed into a contact structure within the space of confoliations. The Giroux correspondence refines this: 

\begin{proposition} \label{GirouxCorresp}
(The Giroux Correspondence\cite{Giroux2002})  For $M^3$, closed and oriented, there is a 1-1 correspondence  between:
{oriented contact structures up to isotopy} and {open book decompositions up to positive stabilization}. 
\end{proposition}

\begin{remark} \label{linFFinBindingNhbd}
In the context of the toy inverse problem, it is important to note that as one approaches a current-carrying wire whose radius of curvature is huge compared to its radius, the magnetic field diverges like $d^{-1}$ where $d$ is the distance to the the center of the wire. This shows that a wire with vanishing radius has an infinite amount of magnetic energy per unit length. Hence, it will be imperative to consider our toy problem and the use of the Giroux correspondence in the context  of finite energy constraints. This is easily handled by replacing the binding of the open book by a tubular neighborhood and solving an eigenvalue problem for a self-adjoint curl operator on the tubular neighborhood. The Lagrangian boundary conditions used to render the curl operator self-adjoint are the ones that correspond to having the trace of the magnetic field correspond to a closed 1-form when pulled back to the boundary. Such boundary conditions are unique up to a choice of Lagrangian subspace of the first homology group of the boundary of the tubular neighborhood when considered as a symplectic vector space\cite{HKT}. The key point to remember is that this choice of self-adjoint boundary condition is tied to the ``positive stabilization" in the Giroux correspondence.
\end{remark}

 \begin{theorem}\label{InducedCorrespondence}
In the absence of ``reconnection points" (as in remark(\ref{PerfectEqNoReconnection})), the Giroux correspondence of  proposition(\ref{GirouxCorresp}) sets up a correspondence between:
\begin{itemize}
\item Solutions to the toy inverse problem as characterized by corollary(\ref{ReformulatedToyProbSolution}) and
\item Near force-free magnetic fields as characterized by Lemma(\ref{near f-f n confoliations}).
\end{itemize}
For fixed magnetic energy, by remarks (\ref{HelicityEVP},\ref{BigEVPplans}, and \ref{linFFinBindingNhbd})  the ``positive stabilization" of the Giroux correspondence articulates the twisting of  a linear force-free magnetic field on a fixed tubular neighborhood of an open book's  binding.  
\end{theorem}

\section{Interlude: The rhetorical quest for ``intuitive cuts"}

Consider the distinction between {\em intuitive but not computable}, and {\em computable but unintuitive}. Let $G$ be the fundamental group of a compact orientable four-manifold. It is a finitely presented group whose generators and relations have a simple interpretation in terms of surgery, yet  basic questions about it  are not decidable in the sense of Turing. This connection is both the basis for proving that there can be no algorithm to classify 4-manifolds up to diffeomorphism, and a prime example of an algebraic structure which is {\em intuitive but not computable}. On the other hand, homology, in any dimension, is computable via the Smith normal form, but as Pontryagin and Thom showed, it deals with cycles which are not necessarily interpretable in terms of manifolds. Homology is a prime example of something which is {\em computable but unintuitive}. In 3-d Poincar\'e duality demands that cuts need to be oriented cycles but, in principle, little else is required. These orientable cycles are computable, but they are not intuitive (especially as subcomplexes of a 3D finite element mesh!) Realizing cuts as orientable embedded manifolds by exploiting the properties of maps into circles doesn't seem to be rooted in intuition and the underlying mathematics can seem divorced from applications, yet history shows otherwise. Hence, we are challenged to understand why certain effective tools for problem solving are dismissed as unintuitive. Taking up this challenge will lead us to valuable insights concerning the work of late-onset blind topologists. Consider ... \\

{\bf A Rhetorical Question: What are nice cuts?}\\

Seeking an answer to this rhetorical question is a game of specifying more constraints on a cut which would still make cuts easy to compute but also intuitive-  especially if this can be achieved with a minimal amount of computational overhead. So, consider three possible suggestions for what a nice cut can be:
\begin{enumerate}
\item Compact Orientable Embedded Manifolds with Boundary (COEMBs)?
\item Genus minimizing COEMBs? 
\item Minimal area COEMBs? 
\end{enumerate}
Some valuable insights follow from preliminary observations:

$\bullet$ The disc is a cut for the textbook example of a circular planar loop and it enjoys all of these properties. However, textbook examples are of little use for probing the rhetorical question. The algorithm based on a harmonic map into the circle, proposed and implemented decades ago (\cite{cuts4}, \cite{gross}), ensures that the cut is a COEMB and finds it in polynomial time with polynomial order quadratic in the worst case (so far so good). Realizing cuts as embedded manifolds is not a logical necessity and so a pragmatic person may question the need  to have embeddings in the first place. However embedded level sets of an elliptic equation come as a by-product of the algorithm and are relatively easy to visualize and understand� (unlike immersions). Unfortunately, there are no conjectures about having a level set which is either genus or area minimizing.

$\bullet$ Given the discussion of the topological characterization of force-free magnetic fields (via Thurston and Alexander norms), it would be wonderful if there were a robust algorithm to produce a cut realized as a genus minimizing COEMB. Unfortunately, Thurston et.al.\cite{HLT} have shown that there is no polynomial time algorithm that produces a COEMB of minimal genus.   In general we can replace the harmonic map functional by the $L^{p}$ norm of the magnitude of the gradient. When we do this for the $L^{1}$ norm, the solution gives an area minimizing current (in the sense of de Rham). However this does not ensure that the solution is unique or an area minimizing COEMB and the discretized problem leads to a linear programming problem for which no worst-case polynomial time algorithm exists. 

$\bullet$ Thurston and Almgren \cite{AlmThurst}, produced unknotted curves which bound only surfaces of high genus within their convex hulls. Since any minimal surface lies in the convex hull of its boundary, this shows that in general one cannot expect to find a COEMB which is both area minimizing and minimal genus. Incidentally, their examples involve unknotted smooth curves with relatively small curvature!

From these deliberations on what a nice cut might be, it is clear that in developing a general algorithm, the harmonic map approach not only side-steps the pitfalls of exponential time algorithms but the solution time is quadratic in the worst case, and practical in practice. The nonexistence of polynomial time algorithms to handle area or genus minimization reinforced the  ``computable but unintuitive" nature of this solution, yet those who sought ``a nicer  and more intuitive" solution still grieve! 

\begin{remark} \label{four late-onset blind}
The work on cuts, as presented so far, exposed us to the work of Pontryagin and Giroux. As we wrestle with the more intuitive questions, we come face to face with the work of other late-onset blind geometers and topologists.  First, the problem of area minimization leads to Plateau's problem. Joseph Plateau was a Belgian solar physicist who was completely blind by the age of forty as a result of staring at the sun. He lived to the age of eighty and his work on capillary forces and minimal surfaces was undertaken long after he was blind! Second, in the process of understanding embedded and immersed cuts the problem of sphere eversion emerges\cite{FrancisMorin79}. That is, turning the sphere inside out homotopically through immersions.  It was Bernard Morin who first showed the world how to visualize the sphere eversion. What was strange was that he was blind from the age of six. 
\end{remark} 

Bottom line: In the quest for {\em nice cuts}, the above paragraphs drops the names of four late-onset blind  topologists. A conclusion one can draw from this exploration is that, in the process of coming up with efficient and effective algorithms to compute cuts, one discovers key topologists and geometers who gave us the tools to construct and visualize cuts and force-free magnetic fields (Pontryagin, Giroux, Plateau, Morin), who were all blinded after the age of five and did their significant mathematical work years after they were blinded. 

\section{Mathematics evolves while embracing Faustian bargains}

In (\cite{Atiyah2002}), Sir Michael Atiyah identified a type of mathematical Faustian bargain in terms of giving up geometry for algebra:
\begin{quote}
 Algebra is the offer made by the devil to the mathematician. The devil says: ``I will give you this
powerful machine, it will answer any question you like. All you need to do is give me
your soul: give up geometry and you will have this marvelous machine." (Nowadays
you can think of it as a computer!) Of course we like to have things both ways; we
would probably cheat on the devil, pretend we are selling our soul, and not give it
away. Nevertheless, the danger to our soul is there, because when you pass over into
algebraic calculation, essentially you stop thinking; you stop thinking geometrically,
you stop thinking about the meaning.
\end{quote}

Although easy to accept intellectually, this became significant for the author when he corresponded with with late-onset blind mathematicians. In this new context the Faustian bargain became one of giving up three dimensional spatial intuition for the convenience of putting arguments on the page, and calculation is only one aspect of the new convenience. The author was surprised to find that Euclidean geometry was not an inspiring subject that appealed to those who would be late-onset blind topologists. Both Morin and Giroux confirmed that it was ``math for lawyers" in that there were two columns, one with pictures and the statement of the the problem, while one had to write out the steps of a proof in the other column. Although important training for writing out arguments and appreciating proofs, it did little to kindle the spark of spatial geometry for these blind students\footnote{The author has to confess that in high school, his spark for three-dimensional thinking was not kindled in geometry class, but in drafting class by exploded views of machine parts, drawings of sections and ``developments". But this leads us to Monge and descriptive geometry}. 

When  Morin and Giroux were pressed on where they felt they had an edge over classmates in geometry class, they came back with the same answer\footnote{Personal communication.}: conic sections! -The teacher would come to class with a wooden cone and ask the students to imagine how circles, ellipses, parabolas and hyperbolas arose as one intersected the cone with a plane. The students would squirm and complain that the teacher should stop making the subject complicated and just draw a picture on the blackboard, while the blind students would wonder what all the fuss was about since it was all obvious! This brings us to ...\\

 The topological Faustian bargain:\\ 
{\bf Give up genuine spatial intuition for the convenience of the page.} \\

It is clear that Giroux and Morin both avoided this topological Faustian bargain. Not by choice, but destiny. In knot theory, they avoid knot projections, consider knots intrinsically; sometimes as the boundary of their Seifert surfaces. To them knot projections are for those who compute knot invariants without soaking in the three dimensional context\footnote{Not all projections are being trivialized here. Projecting four dimensions into three  is natural if one wants an obvious proof that curves can be ``unknotted" in four dimensions. Such a projection is in the style of the late-onset visually impaired.}.

\begin{remark} \label{Bott Smale Morin}
Consider the story of sphere eversion in the context of this Faustian bargain. Stephen Smale, a thesis student of Raoul Bott, proved the result\cite{Smale58}. At the time of Smale's thesis Raoul Bott was an electrical engineer who had metamorphosed into an iconic topologist, known for  Bott periodicity via the application of Morse theory to the loop spaces of the classical Lie groups. So, it is startling that Bott refused to {\em believe} Smale's proof {\em because he couldn't see it!} Not being able to find a gap in the proof, Bott sought strategies to prove the opposite result, but to no avail.  Finally, through the clay models and careful arguments of Bernard Morin, Bott was able to see {\em and feel} the result\footnote{See Silvio Levy (http://www.geom.uiuc.edu/docs/outreach/oi/history.html) for history.}.

The moral is: The Bott ignored Smale's legitimate proof of sphere eversion and exhibited uncharacteristically stubborn and simple-minded behavior when he expected the problem to be resolved by the topological Faustian bargain and it wasn't. On the other hand, Morin did not have the any perceived advantage in succumbing to the temptation offered by the page\footnote{ The story doesn't end here. Over the decades many sphere eversion animations have been made, but they are not particularly effective for visualizing Smale's proof. It is amazing is that Morin can articulate why each one is ineffective! (He sits with people who find the movies ineffective and asks questions as they repeatedly watch the movies. As he teaches them the material, he gets a feel for what was missing in the movie!)}.
\end{remark}

\subsection{Three groups who have avoided the topological Faustian bargain}

It is difficult to talk about this Faustian bargain since written communication and planar diagrams had been ingrained in our culture since well before the Enlightenment. Mathematicians prove a new result by any means available and it is unthinkable that a mathematician would give up calculation, written communication and 2-d visualization for no good reason. Productive mathematicians also find it a distraction to ask how they prove  results given their nationality, gender, physiology, or on account of some other characteristic such as a disability. If one tries to identify those who avoid the topological Faustian bargain, one must start with groups of mathematicians whose use of spatial reasoning is not inspired by the page to begin with, and identify topics that are not amenable to representation on a page. The author speculates there are at least three such groups.

The first group consists of {\em Astronomers}. Although astronomy is highly automated today with many specialists, in previous centuries it was not only at the cutting edge of theoretical physics and experimental technique, but also on the cutting edge of the mathematics. Furthermore, those on the forefront of all these fields were also the ones peering through the telescopes! A planet's trajectory is a curve in a four-dimensional space-time. A star-gazing astronomer tries to make sense of the cosmos for which there is no simple ``model reduction" yet accurate solutions to ill-posed problems are required if one wants to keep track of  individual bodies. The following three astronomers are known primarily for their mathematical work:

$\bullet$ {\em Johan Friedrich Pfaff}. His 1786 Dissertation had to do with celestial mechanics. Pfaff is well known mathematically for his 1812 monograph on ``Pfaffian" (or differential) forms. He figures prominently in the theory of integrability (Frobenius' theorem) and the history of the generalized Stokes' theorem.  Two of Pfaff's students are Gauss and Moebius.

$\bullet$ {\em Carl Friedrich Gauss.} His 1799 Dissertation had to do more with topology than astronomy, but astronomy, geodesy and terrestrial magnetism were the focus of Gauss' ``day job". At least four of Gauss' students are trained as astronomers (Bessel, Gerling, Encke, and von Staudt),  and they generate almost all of Gauss' academic descendants. Gauss estimated the trajectory of Ceres and developed many powerful mathematical techniques in the analysis of astronomical data. Gauss'  mentoring of geometers and topologists (e. g. Listing, Wichmann and Riemann) happened late in his career or in ``retirement".
 
$\bullet$ {\em August Ferdinand Moebius}. His 1815 Dissertation was on astronomy but the theses of his students were uniformly on geometry and astronomy. A moon crater and an asteroid are named after him. In retirement, he and J. B. Listing discovered the Moebius band independently.  
Although he has major contributions to number theory, his contributions to geometry include introduction of homogeneous coordinate into projective geometry, and the Moebius transformation.

Two dimensions may suffice for a textbook explanation of Kepler's contributions, but subsequent developments in astronomy required much more. These astronomers had a great impact on the techniques of n-dimensional geometry\footnote{Their significant results are typically found in diaries and correspondence since Crelle's Journal, the first significant mathematics journal not tied to an Academy, was founded in 1826.}. When one looks at topological aspects of Gauss' work (e.g. fundamental theorem of algebra via winding numbers, the Gauss-Bonnet theorem, and linking numbers) it is clear that he avoided the topological Faustian bargain but difficult to see how since his results often appear in final form as if on stone tablets after a trip up some mountain!

The second group consists of $20^{th}$ century {\em rock climbers}. Although {\em rock climbing topologist} conjures up many names, let's focus on three: James Waddell Alexander II (1888-1971)\footnote{ Alexander's Chimney in Rocky Mountain National Park (CO, USA) is named after him.}, Georges de Rham (1903-1990) and Hassler Whitney (1907-1989). We already encountered mathematical results of the first two, ``Whitney forms" are well-known to those in computational electromagnetics, but this doesn't scratch the surface of the mathematical contributions of each. Each also made significant contributions to rock climbing\footnote{For anecdotes about the interaction of the three climbing together, see page 13 of\cite{WhitneyI}.}. The pressing question is: What does rock climbing have to do with topology? It cannot be answered here but consider the intense planning and decision making required in rock climbing to get a hint. Planning a climb involves local and global decisions made both in advance and on the spot. Path planning is crucial since descending is usually far more dangerous than ascending. Decisions are made on the fly based on both visual and tactile information; getting a foothold around a corner while feet are obscured and the blood circulation in one's hands is rapidly decreasing. Finally, life and death issues revolve around a wide variety of tiny gadgets and tying the correct knot at the right time!

The third group consists of {\em late-onset visually impaired topologists}. It is only in the third group that one can make a reasonable hypothesis about neuroplasticity. 

\subsection{Early vs. late-onset blind mathematicians;  neuroplasticity}

Given the topological characterization of near force-free magnetic fields, the rhetorical search for ``nice cuts"  and the algorithms to compute them, the topological Faustian bargain and those who may have evaded it, �we now hope to avoid aspects of  intuition and conventional wisdom when considering a difficult subject. One cannot systematically survey blind mathematicians, since it would be impossible to find a statistically significant sample within a group of brilliant topologists. However, we are equally blind in higher dimensions\footnote{ Notions like ``transversality" from differential topology, when applied to higher dimensional spaces, give us a feel for how one can do mathematics in a context where sighted and nonsighted are equally blind.} and there are two things one can do:
\begin{enumerate}
\item List famous blind mathematicians to identify characteristics of topologists.
\item Correspond with blind topologists, hoping to avoid misconceptions.
\end{enumerate}
Eight blind mathematicians are listed below\footnote {When discussing this with David Mumford years ago, he suggested that  notoriously short sighted topologists who seldom proof read their papers should be included on this list. In particular he was thinking of Henri Poincar\'e, Bill Thurston and possibly Leonhard Euler. Indeed considering the evolution of three-dimensional topology, the names of Poincar\'e and Thurston figure prominently. However, as we shall soon see, they are different from the late-onset blind (Plateau, Pontryagin, Morin and Giroux) in that their disabilities do not force them make a clean break from the Faustian bargain as presented by Atyiah. Furthermore, one cannot reasonably speculate on the role of neuroplasticity in their work since they they had their eyesight at the beginning of their professional careers and they lack a distinct event that contributed to a total loss of eyesight.}. Although statistically insignificant, it exhibits excellent correlation between onset of blindness and field of research:

$\bullet$ Nicholas Saunderson (1682-1739), Cambridge University's $4^{th}$ Lucasian Prof. (Newton was the $2^{nd}$ and a colleague), perhaps the first to discover Bayes' Theorem. Performed spectacular calculations on an abacus he developed. Blind from birth.

$\bullet$ Joseph Antoine Ferdinand Plateau (1801-1883), see remark (\ref{four late-onset blind}), Belgian physicist (who also studied math.), completely blind by age 42, Contributions to minimal surfaces (Plateau's problem), and capillary action were initiated later on in life.

$\bullet$ Louis Antoine (1888-1971) was blinded at age 29 in WWI. After the war Henri Lebesgue, his thesis advisor directed him toward  topology since it was suited to those who do not like to read. He published on {\em Antoine's Necklace}� in 1921. A high school math teacher by training, he went to grad school and taught in a university only after losing his sight.

$\bullet$ Lev Semenovitch Pontryagin (1908-1988); many contributions to topology but also recognized for work in applied mathematics done after 1952. Blinded by an explosion at age 14, he started a serious study of topology at age 19. Known for Pontryagin duality, solving Hilbert's fifth problem for abelian groups, the Pontryagin-Thom construction (see remark (\ref{Pontryagin})), using Morse theory to compute the Poincar\'e polynomials of classical compact Lie groups, and Pontryagin (characteristic) classes.
 
$\bullet$ Bernard Morin (1931-), see remark (\ref{Bott Smale Morin}), blinded at age six. Known for work in low dimensional topology, sphere eversions and Boy's surface parameterizations.

$\bullet$ Abraham Nemeth (1918-2013) Math. Prof., known primarily for mathematical extensions to the Braille code and the Math Speak Initiative. Blind from Birth.

$\bullet$ Lawrence Baggett, blinded at age five, known for work in harmonic analysis.

$\bullet$ Emmanuel Giroux, blind by age eleven, contact geometer, known for {\em the Giroux correspondence} relating {\em open books} to {\em isotopy classes of contact structures}.

\bigskip
These eight mathematicians fit nicely into three groups:
\begin{enumerate}
\item {\bf Saunderson, Nemeth}. Blind from birth, their mathematics is mainstream, and doesn't have a distinctly spatial aspect.
\item {\bf Pontryagin, Plateau, Morin and Giroux} are unique; blinded after the age of five, they learned their type mathematics more than five years after losing their sight. The mathematics they are known for is not computational or conjectural,  but involves concrete geometric or topological constructions.
\item {\bf Antoine and Baggett} are analysts. One is tempted to put Louis Antoine in the second group since he was blinded after age five but he differs from those in that group since his main research results were obtained within five years of being blinded. Although his work is considered topology, his techniques of proof are more akin to those of an analyst than a topologist \footnote{The author is grateful to Bernard Morin for this insight}. Being blinded at the borderline age of the late-onset blind, one is also tempted to put Lawrence Baggett in the second group but he mastered the tools of a conventional mathematician, and made his mark as an analyst.  
\end{enumerate}

Of these three groups, the second is distinct from the first and third by the crucial role played by cross-modal plasticity! (Cross-modal plasticity is a specific example of neuroplasticity and it is particularly significant in the present context). One can see that the mathematicians in the second group spent well over 10,000 hours in a world where they integrated their visual systems with their other senses as they did the things that kids do, and then another 10,000 hours in a world where {\em path planning without visual input} was key, before obtaining major mathematical results. In fewer words we might advance the following ...

\bigskip
{\center{\bf Hypothesis: Late-onset blind mathematicians with 10,000 hours\\ of visually-deprived path planning are topologists.} }

\bigskip
It is impossible to test this hypotheses with a statistically significant sample. If progress in mathematics requires studying the masters and not their students, it behooves us appreciate the consequences of this hypothesis through a study of scarce but useful literature and of living pioneers (without distracting them from their research!) Alternatively, one might find a statistically significant sample by studying situations where sighted people are distracted by their visual systems, and choose to ignore them\footnote{Three such instances come to mind. First I am grateful to Douglas Hofstader for pointing me to chess masters who prefer to keep track of the game by imagining the board rather than by looking at it (blind chess). Second, one could point to the practice of studying a map before embarking on a trip to a new destination, and not looking at it while traveling. Third, there are legends of geometers, most notably Jacob Steiner, who successfully taught geometry in the dark.}. To see why neuroplasticity should play a central role in the hypothesis, consider more of Sir Michael Atiyah's wisdom (\cite{Atiyah2002}): 
\begin{quotation} 
Our brains have been constructed in such a way that they are extremely concerned
with vision. Vision, I understand from friends who work in neurophysiology, uses
up something like 80 or 90 percent of the cortex of the brain. There are about 17
different centres in the brain, each of which is specialised in a different part of the
process of vision: some parts are concerned with vertical, some parts with horizontal,
some parts with colour, or perspective, and finally some parts are concerned with
meaning and interpretation. Understanding, and making sense of, the world that we
see is a very important part of our evolution. Therefore, spatial intuition or spatial
perception is an enormously powerful tool, and that is why geometry is actually such
a powerful part of mathematics|not only for things that are obviously geometrical,
but even for things that are not. We try to put them into geometrical form because
that enables us to use our intuition.
\end{quotation} 
Late onset visual impairment and subsequent cross-modal plasticity is key to the hypothesis because the late onset blind had the neural connections of a typical visual system which were coopted for a world of careful nonvisual task planning\footnote{{\em Path planning without visual input}, suggests an analogy: Walking with visual feedback is like driving absorbed and enslaved by a GPS; a person walking effectively without visual feedback is like a person who looks at a map once and drives to an unfamiliar place without a GPS.}.

\section{Conclusions and outlook}
The main result of the paper is theorem (\ref{InducedCorrespondence}). It exploits the Giroux correspondence to relate two characterizations of near force-free magnetic fields (in the absence of reconnection points). One characterization is a consequence of the definition of  ``near force-free magnetic field", and the other comes from the consideration of a ``toy" inverse problem. Along the way, we noted how  leveraging results from higher dimensional topology isolates abstractions, can lead to cleaner development of algorithms, and that these tools often have their roots in three dimensions- but in the hands of late-onset blind topologists subject to an additional hypotheses pertaining to neuroplasticity. This additional hypothesis is carefully developed in the context of a ``topological Faustian bargain".

This research can be taken in many directions. The function spaces arising from equivalence relations arising in topology yield natural roles for the calculus of variations. Historically, this is the story of Morse theory and Hodge theory. For cuts it was harmonic maps, for sphere eversions, it can be Willmore energy, for linear force-free magnetic fields it is an eigenvalue problem for a helicity functional. All these tie the finite element method to challenging problems in visualization.

The first side lesson  is: heuristics and intuition invite lousy algorithms and denial about the applicability of mathematics. Pedagogically, there's the issue of programming environments for late-onset blind STEM students, the use of Latex, and the relative merits of looking at pictures versus creating one's own pictures and models. 

\section*{Acknowledgments} The author thanks  Emmanuel Giroux, Bernard Morin,  Lorraine Norwich and Tomasso Toffoli who listened to his speculations and responded with solid insights. Also, Alberto Valli for fostering unconventional applied math,  Brendan Mattingly for staying ``on the same page as the class" by learning Latex (and by his example getting the class to reciprocate), and Thomas Garrity for sharing the Latex files for his text with a blind student and taking a genuine interest.

\medskip
Received xxxx 20xx; revised xxxx 20xx.
\medskip

\end{document}